# STABILITY RESULTS FOR SIMPLE TRAFFIC MODELS UNDER PI-REGULATOR CONTROL

Iasson Karafyllis[*] and Markos Papageorgiou[**]

[*]Dept. of Mathematics, National Technical University of Athens,
Zografou Campus, 15780, Athens, Greece (email: iasonkar@central.ntua.gr )

[**]School of Production Engineering and Management, Technical University of Crete,
Chania, 73100, Greece (email: markos@dssl.tuc.gr )

**Abstract**

This paper provides necessary conditions and sufficient conditions for the (global) Input-to-State Stability property of simple uncertain vehicular-traffic network models under the effect of a PI-regulator. Local stability properties for vehicular-traffic networks under the effect of PI-regulator control are studied as well: the region of attraction of a locally exponentially stable equilibrium point is estimated by means of Lyapunov functions. All obtained results are illustrated by means of simple examples.

## 1. Introduction

There are a number of relatively simple controlled processes within vehicular-traffic networks or water networks, which share the following characteristics:

- The kernel of the process is some sort of "reservoir" (e.g. an urban road network, a freeway stretch, a water reservoir or basin) which accumulates inflows and outflows; the reservoir features a limited storage capacity.
- There is a controllable but constrained inflow; the inflow may be released at some distance from the reservoir, in which case it reaches the reservoir with a corresponding time-delay.
- There may be additional uncontrollable inflows.
- The outflow depends on the reservoir storage in a nonlinear way; there may be some modeling uncertainty in the related function.
- The control goal is to operate the system near a pre-specified storage level.

Examples of such controlled processes include local freeway ramp metering [13], gating control of urban network parts [9], merging traffic control [14], variable speed limit control on freeways [4], water level and water flow control [11,12]. In some cases, these elementary systems may be interconnected to form bigger composite systems, as, e.g., in the cases of multiple urban network parts [1] or irrigation networks [3].

The mentioned characteristics indicate that these elementary processes may be modeled as discrete-time time-delayed constrained nonlinear first-order systems. A PI-type regulator is usually employed for system control in practice; whereby the regulator parameters are selected after model linearization around the desired set-value, using classical linear sample-data concepts. It should be noted that, in the case of traffic systems, the nonlinear function connecting the outflow with the reservoir storage is typically a concave uni-modal function featuring a maximum, which usually corresponds to the desired operation state.

Although these systems are usually operating reasonably well in practice, it is interesting to have a second look at them from a nonlinear analysis point of view. Specifically, we are interested in deriving local and global stability results for the PI-controlled nonlinear models, which is the main scope of this paper (albeit without consideration of possible input delays). Eventually, we are interested in deriving nonlinear stabilizers and, finally, in considering control of bigger composite systems.

Consider the following 1-dimensional discrete-time control system, which is representative for all examples of elementary practical systems mentioned earlier (without input delay):

$$x^+ = x - f(d,x) + \min(u+v, a-x)$$
$$x \in [0,a], u \in [b_{\min}, b_{\max}], d \in D, v \in \Re_+ \quad (1.1)$$

where $D \subseteq \Re$ is a compact set, $a > 0$, $0 \leq b_{\min} < b_{\max}$ are constants and

$$f \in C^0(D \times [0,a]; \Re) \text{ with } 0 \leq f(d,x) \leq x \text{ for all } (d,x) \in D \times [0,a]. \quad (1.2)$$

System (1.1) describes the time evolution of a traffic (or water) system, where $x \in [0,a]$ is the current number of vehicles (storage) in the network, $f(d,x)$ is the (uncertain) outflow function, $a > 0$ is the capacity of the network and $v \in \Re_+$ is the input that reflects the uncontrollable inflow. System (1.1) under assumption (1.2) is a well-defined control system which satisfies $x^+ \in [0,a]$ for all $(x,u,d,v) \in [0,a] \times [b_{\min}, b_{\max}] \times D \times \Re_+$.

In order to state the control problem, we assume that:

**(H1)** *There exists* $(x^*, u^*, v^*) \in (0,a) \times (b_{\min}, b_{\max}) \times \Re_+$ *such that* $f(d, x^*) = u^* + v^* < a - x^*$ *for all* $d \in D$.

In other words, we assume that $x^* \in (0,a)$ is an equilibrium point for system (1.1) with $v \equiv v^*$ and $u \equiv u^*$.

The PI regulator is the dynamic feedback law that is given by the equation:

$$u(t) = \max\left(b_{\min}, \min\left(b_{\max}, u(t-1) - k_1(x(t) - x(t-1)) - k_2(x(t) - x^*)\right)\right) \quad (1.3)$$

where $k_1, k_2$ are constants. The closed-loop system (1.1) with (1.3) is described by the 3-dimensional discrete-time system:

$$x^+ = x - f(d,x) + \min(u+v, a-x)$$
$$y^+ = x$$
$$w^+ = u \quad (1.4)$$
$$u = \max\left(b_{\min}, \min\left(w - k_1(x-y) - k_2(x-x^*), b_{\max}\right)\right)$$

with state space $S = [0,a]^2 \times \Re$ (i.e., $(x,y,w) \in [0,a] \times [0,a] \times \Re$) and inputs $(v,d) \in \Re_+ \times D$. The point $(x^*, x^*, u^*) \in (0,a)^2 \times (b_{\min}, b_{\max})$ is an equilibrium point of system (1.4). In this work, we answer the following questions concerning the PI regulator:

1) What are the conditions that guarantee the (global) Input-to-State Stability (ISS) property with respect to the external input $v \in \Re_+$ uniformly in $d \in D$ for system (1.4)?
2) What are the conditions that guarantee local exponential stability for the equilibrium point $(x^*, x^*, u^*) \in (0,a)^2 \times (b_{\min}, b_{\max})$ in the disturbance-free case, i.e., when $f(d,x)$ is independent of $d \in D$ and $v \equiv v^*$?
3) What is the region of attraction when the equilibrium point $(x^*, x^*, u^*) \in (0,a)^2 \times (b_{\min}, b_{\max})$ is locally exponentially stable in the disturbance-free case?

As expected, the answers to the above questions are related. The notion of ISS for discrete-time systems was studied in [7] and this notion is adopted here, although the system that we study (namely system (1.4)) evolves in a restricted state space (in $S = [0,a]^2 \times \Re$) and not in $\Re^3$ (all the results in [7] are for discrete-time systems evolving in $\Re^n$). When the ISS property is applied to system (1.4) with $v \equiv v^*$, then it becomes identical to the notion of the Robust Global Asymptotic Stability (see [8]). Stability properties for discrete-time systems with restricted state spaces are studied in [17].

The structure of the paper is as follows: Section 2 is devoted to the answer of questions (2) and (3) above, while Section 3 addresses question (1) above. All obtained results are illustrated by means of some simple examples in Section 4. Finally, the concluding remarks of the present work are provided in Section 5.

**Notation.** Throughout this paper, we adopt the following notation:

* $\Re_+ := [0,+\infty)$.
* By $C^0(A;\Omega)$, we denote the class of continuous functions on $A \subseteq \Re^n$, which take values in $\Omega \subseteq \Re^m$. By $C^k(A;\Omega)$, where $k \geq 1$ is an integer, we denote the class of functions on $A \subseteq \Re^n$ with continuous derivatives of order $k$, which take values in $\Omega \subseteq \Re^m$.

## 2. Local Results

This section is devoted to the analysis of local exponential stability for the disturbance-free version of the closed-loop system (1.1) with (1.3). The local stability analysis of the disturbance-free version of the closed-loop system (1.1) with (1.3) is equivalent to the stability analysis of the following system:

$$\begin{aligned} x^+ &= x - f(x) + \min(w + v^*, a - x) \\ w^+ &= P\big(w + \sigma f(x) - \sigma \min(w + v^*, a - x) - k_2(x - x^*)\big) \\ (x,w) &\in [0,a] \times \Re \end{aligned} \quad (2.1)$$

where $f \in C^1([0,a];\Re)$ satisfies (H1), $\sigma = k_1 + k_2$ and

$$P(x) = \max\big(b_{\min}, \min(x, b_{\max})\big), \text{ for all } x \in \Re \quad (2,2)$$

Indeed, it should be noticed that the solution $(x(t), y(t), w(t))$ of (1.4) satisfies equations (2.1) for $t \geq 1$, since the equations

$$\begin{cases} y(t+1) = y(t) - f(y(t)) + \min(w(t) + v^*, a - y(t)) \\ w(t+1) = P\big(w(t) + \sigma f(y(t)) - \sigma \min(w(t) + v^*, a - y(t)) - k_2(y(t) - x^*)\big) \end{cases}$$

hold for $t \geq 1$. The equilibrium point $(x^*, u^*)$ of (2.1) is in the interior of the region $\big\{(x,w) \in [0,a] \times \Re : w + x \leq a - v^*, b_{\min} \leq (1-\sigma)w + \sigma f(x) - \sigma v^* - k_2(x - x^*) \leq b_{\max}\big\}$ (recall (H1); notice that the right hand side of (2.1) is continuously differentiable on the interior of the previously mentioned region), and therefore it follows that the equilibrium point $(x^*, u^*)$ of (2.1) is locally exponentially stable if and only if all roots of the equation

$$s^2 - (2 - f'(x^*) - \sigma)s + (1 - f'(x^*) - \sigma + k_2) = 0 \quad (2.3)$$

are strictly inside the unit ball (see Chapter 5 in [16], Chapter 4 in [10] and the necessary extensions to the case of local exponential stability). In other words, one of the following conditions is equivalent to local exponential stability of the equilibrium point $(x^*, u^*)$ of (2.1):

**(I)** $(2 - f'(x^*) - \sigma)^2 < 4(1 - f'(x^*) - \sigma + k_2) < 4$ (complex roots)

**(II)** $\big|2 - f'(x^*) - \sigma\big| < 2$ and $\big|2 - f'(x^*) - \sigma\big| - 1 < 1 - f'(x^*) - \sigma + k_2 \leq \dfrac{(2 - f'(x^*) - \sigma)^2}{4}$ (real roots).

In order to give an estimation of the region of attraction of the equilibrium point $(x^*, u^*)$ of (2.1), we need to perform a Lyapunov analysis. The Lyapunov function for the case of a two dimensional linear discrete-time system with characteristic polynomial $s^2 + bs + c = 0$ and real roots strictly inside the unit ball is the function $V = |\xi| + M\left|\xi^+ + \frac{b}{2}\xi\right|$ for an appropriate constant $M > 1$, where $\xi$ is the deviation of any state of the system from its equilibrium value. Therefore, the following function

$$V(x,w) := |x - x^*| + M|w + v^* - f(x) + (1-g)(x - x^*)|, \text{ for all } (x,w) \in [0,a] \times \Re \quad (2.4)$$

where $M > 1$ and $g \in \Re$ with $|g| < 1$ are constants to be determined, is a candidate Lyapunov function for (2.1). An estimation of the region of attraction $A \subseteq [0,a] \times \Re$ of the equilibrium point $(x^*, u^*)$ of (2.1) is the sublevel set

$$\Omega_\rho = \{(x,w) \in [0,a] \times \Re : V(x,w) < \rho\} \subseteq A \quad (2.5)$$

where $\rho > 0$ is a constant, for which the inequality $V(x^+, w^+) < V(x,w)$ holds for all $(x,w) \in \Omega_\rho$ (see Chapter 4 in [10] and [5]).

The following proposition gives an estimation of the region of attraction for certain regions in the parameter space $k_1, k_2$.

**Proposition 2.1:** *Suppose that the equilibrium point $(x^*, u^*)$ of (2.1) is locally exponentially stable. Moreover, suppose that $|k_2 - q^2| < (|1-q|-1)^2$, where $q = \frac{\sigma + f'(x^*)}{2} \in (0,2)$ and $\sigma = k_1 + k_2$. Let $\eta > 0$ be a constant for which $\max\{|f'(x) - f'(x^*)| : x \in [0,a], |x - x^*| \leq \eta\} = L$ and such that $|f'(x) - f'(x^*)| < L$, for all $x \in [0,a]$ with $|x - x^*| < \eta$, where $L := \frac{(|1-q|-1)^2 - |k_2 - q^2|}{|q| + 1 - |1-q|}$. Define*

$$\rho = \min\left(\eta, \frac{\min(b_{\max} - u^*, u^* - b_{\min})}{\max\left(\frac{|1-2q+f'(x^*)|}{M}, L + |k_2 + (1-2q)q + (q-1)f'(x^*)|\right)}, \frac{a - v^* - u^* - x^*}{1 + L + |q - f'(x^*)|}\right) \quad (2.6)$$

*where $M = \frac{|q| + 1 - |1-q|}{|q|(1 - |1-q|) + |k_2 - q^2|}$. Consider the solution $(x(t), w(t)) \in [0,a] \times \Re$ of (2.1) with initial condition $(x(0), w(0)) \in \Omega_\rho = \{(x,w) \in [0,a] \times \Re : |x - x^*| + M|w + v^* - f(x) + q(x - x^*)| < \rho\}$. Then $\lim_{t \to +\infty}(x(t), w(t)) = (x^*, u^*)$.*

Proposition 2.1 does not provide an estimation of the region of attraction for all pairs of values of the parameters $k_1, k_2$, for which the equilibrium point $(x^*, u^*)$ of (2.1) is locally exponentially stable. This is clearly shown in Figure 1 below.

Proposition 2.1 provides a conservative estimation of the region of attraction. In order to obtain a less conservative estimation of the region of attraction, we can also use the following proposition. Proposition 2.2, which delivers a different region of attraction. Specifically, Proposition 2.2 can be applied to values of the parameters $k_1, k_2$, for which Proposition 2.1 can be applied as well, and the overall estimation of the region of attraction corresponds to the union of the regions of attraction resulting from each proposition (see Example 4.2 below).

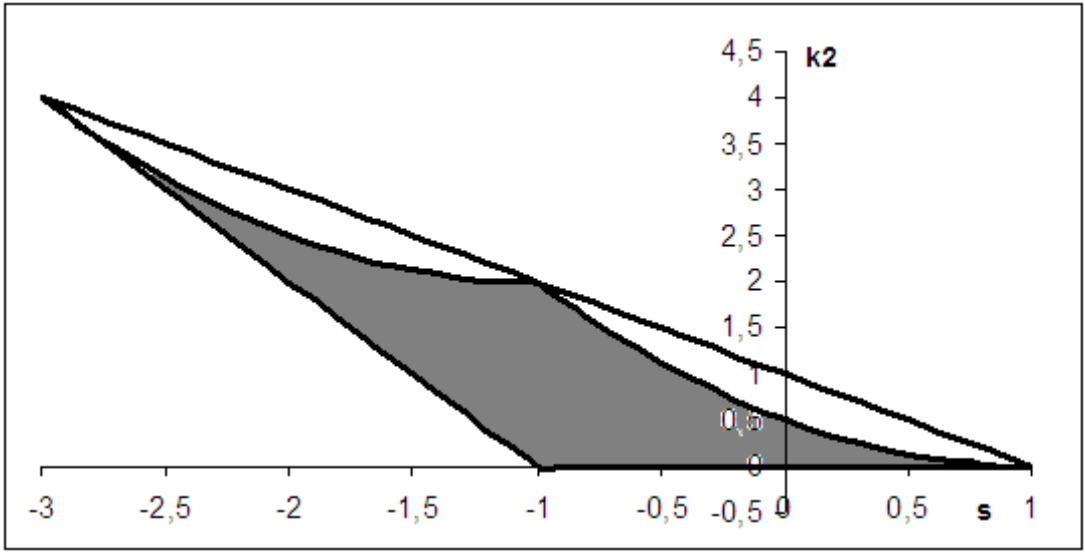

**Fig. 1:** The triangle in the parameter space $s - k_2$, where $s = 1 - f'(x^*) - k_1 - k_2$, for which the equilibrium point $(x^*, u^*)$ of (2.1) is locally exponentially stable. The grey region is the region in the parameter space $s - k_2$ for which Proposition 2.1 can be applied.

**Proposition 2.2:** *Suppose that the equilibrium point $(x^*, u^*)$ of (2.1) is locally exponentially stable. Moreover, suppose that $0 < k_2 < 2$ and $\left| f'(x^*) - 1 + k_1 \right| < \dfrac{1 - |1 - k_2|}{k_2 + 1 - |1 - k_2|}$. Let $\eta > 0$ be a constant for which*

$$\max\left\{ |f'(x) - 1 + k_1| : x \in [0, a], |x - x^*| \leq \eta \right\} = \frac{1 - |1 - k_2|}{k_2 + 1 - |1 - k_2|}$$

*and such that $\left| f'(x) - 1 + k_1 \right| < \dfrac{1 - |1 - k_2|}{k_2 + 1 - |1 - k_2|}$, for all $x \in [0, a]$ with $|x - x^*| < \eta$. Define*

$$\rho := \min\left( \eta, \frac{\min(b_{\max} - u^*, u^* - b_{\min})}{\max\left( \dfrac{k_2|1 - \sigma|}{k_2 + 1 - |1 - k_2|}, \dfrac{1 - |1 - k_2|}{k_2 + 1 - |1 - k_2|} + |(1-\sigma)(k_2 - 1)| \right)}, \frac{a - v^* - u^* - x^*}{1 + \dfrac{1 - |1 - k_2|}{k_2 + 1 - |1 - k_2|} + |\sigma - 1|} \right). \quad (2.7)$$

*Consider the solution $(x(t), w(t)) \in [0, a] \times \Re$ of (2.1) with initial condition $(x(0), w(0)) \in \Omega_\rho = \left\{ (x, w) \in [0, a] \times \Re : |x - x^*| + \dfrac{k_2 + 1 - |1 - k_2|}{k_2} |w + v^* - f(x) + k_2(x - x^*)| < \rho \right\}$. Then*

$$\lim_{t \to +\infty} (x(t), w(t)) = (x^*, u^*).$$

The proofs of Proposition 2.1 and Proposition 2.2 are based on the following facts.

**FACT I:** *Suppose that there exist constants $\eta > 0$, $F \in \Re$, $L > 0$, such that $|f(x) - f(x^*) - F(x - x^*)| \leq L|x - x^*|$, for all $x \in [0, a]$ with $|x - x^*| \leq \eta$. Let $g \in \Re$ and $M \geq 1$ be constants and define*

$$V(x, w) := |x - x^*| + M|w + v^* - f(x) + (1 - g)(x - x^*)|$$

*for all $(x, w) \in [0, a] \times \Re$. If*

$$V(x,w) \le \min\left(\eta, \frac{\min(b_{\max}-u^*, u^*-b_{\min})}{\max\left(\frac{|1-\sigma|}{M}, L+|k_2+(1-\sigma)(1-g)-F|\right)}, \frac{a-v^*-u^*-x^*}{1+L+|1-F-g|}\right) \text{ then } b_{\min} \le (1-\sigma)w - \sigma v^* + \sigma f(x) - k_2(x-x^*) \le b_{\max}$$

and $\min(w+v^*, a-x) = w+v^*$.

**Proof of Fact I:** Indeed, the inequality $V(x,w) \le \min\left(\eta, \frac{\min(b_{\max}+v^*-f(x^*), f(x^*)-v^*-b_{\min})}{\max\left(\frac{|1-\sigma|}{M}, L+|k_2+(1-\sigma)(1-g)-F|\right)}\right)$ implies the inequality:

$$|1-\sigma||w+v^*-f(x)+(1-g)(x-x^*)| + (L+|k_2+(1-\sigma)(1-g)-F|)|x-x^*| \le \min(b_{\max}+v^*-f(x^*), f(x^*)-v^*-b_{\min}).$$

The above inequality in conjunction with the fact that $|f(x)-f(x^*)-F(x-x^*)| \le L|x-x^*|$, for all $x \in [0,a]$ with $|x-x^*| \le \eta$, implies the following inequality:

$$|(1-\sigma)(w+v^*-f(x)+(1-g)(x-x^*)) + f(x)-f(x^*)-F(x-x^*) - (k_2+(1-\sigma)(1-g)-F)(x-x^*)|$$
$$\le \min(b_{\max}+v^*-f(x^*), f(x^*)-v^*-b_{\min})$$

which directly gives:

$$b_{\min}+v^*-f(x^*) \le (1-\sigma)(w+v^*-f(x)+(1-g)(x-x^*)) + f(x)-f(x^*) - (k_2+(1-\sigma)(1-g))(x-x^*) \le b_{\max}+v^*-f(x^*).$$

The above inequality is equivalent to the inequality $b_{\min} \le (1-\sigma)w - \sigma v^* + \sigma f(x) - k_2(x-x^*) \le b_{\max}$.

On the other hand, the inequality $V(x,w) \le \min\left(\eta, \frac{a-v^*-u^*-x^*}{1+L+|1-F-g|}\right)$ implies the inequality:

$$|w+v^*-f(x)+(1-g)(x-x^*)| + (1+L+|1-F-g|)|x-x^*| \le a-v^*-x^*-u^*.$$

The above inequality in conjunction with the fact that $|f(x)-f(x^*)-F(x-x^*)| \le L|x-x^*|$, for all $x \in [0,a]$ with $|x-x^*| \le \eta$, implies the following inequality:

$$|w+v^*-f(x)+(1-g)(x-x^*)| + |x-x^*|$$
$$+|v^*-f(x^*)+f(x^*)-f(x)+F(x-x^*)+(1-F-g)(x-x^*)| \le a-v^*-x^*$$

which directly gives:

$$|w| + |x-x^*| \le a-v^*-x^*.$$

The above inequality implies the inequality $w+x \le a-v^*$, or equivalently, $\min(w+v^*, a-x) = w+v^*$. ◁

**FACT II:** *Suppose that there exist constants $\eta > 0$, $F \in \Re$, $L > 0$ such that $|f(x)-f(y)-F(x-y)| \le L|x-y|$, for all $x,y \in [0,a]$ with $|x-x^*| \le \eta$, $|y-x^*| \le \eta$. Let $g \in \Re$ with $|g| \le 1$ and $M \ge 1$ be constants and define*

$$V(x,w) := |x - x^*| + M|w + v^* - f(x) + (1-g)(x - x^*)| \qquad \text{for all} \qquad (x,w) \in [0,a] \times \Re. \qquad \text{If}$$

$$V(x,w) \leq \min\left(\eta, \frac{\min(b_{\max} - u^*, u^* - b_{\min})}{\max\left(\frac{|1-\sigma|}{M}, L + |k_2 + (1-\sigma)(1-g) - F|\right)}, \frac{a - v^* - u^* - x^*}{1 + L + |1 - F - g|}\right) \qquad \text{then}$$

$$\begin{aligned} V^+ &\leq \left(1 + M|2 - \sigma - g - F| + ML\right)|w + v^* - f(x) + (1-g)(x - x^*)| \\ &+ \left(|g| + M|k_2 + (1 - F - \sigma - g)(1-g)| + ML|1-g|\right)|x - x^*| \end{aligned} \qquad (2.8)$$

where $V^+ = V\left(x - f(x) + \min(w + v^*, a - x), P\left(w - \sigma\min(w + v^*, a - x) + \sigma f(x) - k_2(x - x^*)\right)\right)$.

**Proof of Fact II:** First notice that, by virtue of Fact I, we get $P\left(w + \sigma f(x) - \sigma\min(w + v^*, a - x) - k_2(x - x^*)\right) = (1-\sigma)w - \sigma v^* + \sigma f(x) - k_2(x - x^*)$ and $\min(w + v^*, a - x) = w + v^*$. Using the definition $V(x,w) := |x - x^*| + M|w + v^* - f(x) + (1-g)(x - x^*)|$ and the triangle inequality, we get:

$$\begin{aligned} V^+ &= |x + v^* - f(x) + w - x^*| \\ &+ M|(2 - \sigma - g)w + (\sigma + g - 1)(f(x) - v^*) + (1 - g - k_2)(x - x^*) + v^* - f(x - f(x) + w + v^*)| \\ &\leq |w + v^* - f(x) + (1-g)(x - x^*)| + |g||x - x^*| \\ &+ M|(2 - \sigma - g)z_2 + (1 - g - k_2 - (2 - \sigma - g)(1-g))z_1 + f(x^* + z_1) - f(x^* + gz_1 + z_2)| \end{aligned} \qquad (2.9)$$

for all $(x,w) \in [0,a] \times \Re$ with $V(x,w) \leq \min\left(\eta, \frac{\min(b_{\max} - u^*, u^* - b_{\min})}{\max\left(\frac{|1-\sigma|}{M}, L + |k_2 + (1-\sigma)(1-g) - F|\right)}, \frac{a - v^* - u^* - x^*}{1 + L + |1 - F - g|}\right)$, where $z_1 = x - x^*$, $z_2 = w + v^* - f(x) + (1-g)(x - x^*)$. Since $|g| \leq 1$, $M \geq 1$, it follows from definition $V(x,w) := |x - x^*| + M|w - f(x) + (1-g)(x - x^*)|$:

$$|gz_1 + z_2| \leq |g||z_1| + |z_2| \leq |z_1| + M|z_2| = V(x,w) \leq \eta \text{ and } |z_1| \leq |z_1| + M|z_2| = V(x,w) \leq \eta.$$

Consequently, it follows from the fact that $|f(x) - f(y) - F(x - y)| \leq L|x - y|$, for all $x, y \in [0,a]$ with $|x - x^*| \leq \eta$, $|y - x^*| \leq \eta$:

$$\left|f(x^* + gz_1 + z_2) - f(x^* + z_1) - F(g-1)z_1 - Fz_2\right| \leq L|(g-1)z_1 + z_2|. \qquad (2.10)$$

Inequality (2.8) is a direct consequence of (2.9) and (2.10). ◁

**FACT III:** *Suppose that there exist constants $\eta > 0$, $F \in \Re$, $L > 0$, such that $|f(x) - f(y) - F(x - y)| \leq L|x - y|$, for all $x, y \in [0,a]$ with $|x - x^*| \leq \eta$, $|y - x^*| \leq \eta$. Let $g \in \Re$ with $|g| < 1$ and $M > 1$ be constants with*

$$\frac{1}{1 - L - |2 - \sigma - g - F|} < M < \frac{1 - |g|}{|k_2 + (1 - F - \sigma - g)(1-g)| + L|1 - g|} \qquad \text{and} \qquad \text{define}$$

$$V(x,w) := |x - x^*| + M|w - f(x) + (1-g)(x - x^*)| \quad \text{for all} \quad (x,w) \in [0,a] \times \Re. \quad \text{If}$$

$$V(x,w) \leq \min\left(\eta, \frac{\min(b_{\max} - u^*, u^* - b_{\min})}{\max\left(\frac{|1-\sigma|}{M}, L + |k_2 + (1-\sigma)(1-g) - F|\right)}, \frac{a - v^* - u^* - x^*}{1 + L + |1 - F - g|}\right) \quad \text{then}$$

$$V^+ \leq \lambda V(x,w) \tag{2.11}$$

where $\quad V^+ = V\big(x - f(x) + \min(w + v^*, a - x), P(w - \sigma \min(w + v^*, a - x) + \sigma f(x) - k_2(x - x^*))\big) \quad$ and $\lambda := \max\left(M^{-1} + |2 - \sigma - g - F| + L, |g| + M|k_2 + (1 - F - \sigma - g)(1 - g)| + ML|1 - g|\right) < 1$.

Fact III is a direct consequence of Fact II.

We are now ready to provide the proofs of Proposition 2.1 and Proposition 2.2.

**Proof of Proposition 2.2:** Selecting $F = 1 - k_1$, we notice that the assumption $|f'(x^*) - 1 + k_1| < \frac{1 - |1 - k_2|}{k_2 + 1 - |1 - k_2|}$ guarantees that $|f(x) - f(y) - F(x - y)| \leq \frac{1 - |1 - k_2|}{k_2 + 1 - |1 - k_2|}|x - y|$, for all $x, y \in [0, a]$ with $|x - x^*| \leq \eta$, $|y - x^*| \leq \eta$, where $\eta > 0$ is the constant that satisfies $\max\{|f'(x) - 1 + k_1| : x \in [0, a], |x - x^*| \leq \eta\} = \frac{1 - |1 - k_2|}{k_2 + 1 - |1 - k_2|}$ and $|f'(x) - 1 + k_1| < \frac{1 - |1 - k_2|}{k_2 + 1 - |1 - k_2|}$, for all $x \in [0, a]$ with $|x - x^*| < \eta$.

Let $(x(0), w(0)) \in \Omega_\rho = \left\{(x, w) \in [0, a] \times \Re : |x - x^*| + \frac{k_2 + 1 - |1 - k_2|}{k_2}|w + v^* - f(x) + k_2(x - x^*)| < \rho\right\}$, where $\rho$ is defined by (2.7). Since $f \in C^1([0, a]; \Re)$ and since $|f'(x) - 1 + k_1| < \frac{1}{2}$, for all $x \in [0, a]$ with $|x - x^*| < \eta$, there exists $L \in \left(0, \frac{1 - |1 - k_2|}{k_2 + 1 - |1 - k_2|}\right)$ such that $L := \max\{|f'(x) - 1 + k_1| : x \in [0, a], |x - x^*| \leq \tilde{\eta}\}$, where $\tilde{\eta} = |x(0) - x^*| + \frac{k_2 + 1 - |1 - k_2|}{k_2}|w(0) + v^* - f(x(0)) + k_2(x(0) - x^*)| < \rho$. Consequently, we get $|f(x) - f(y) - F(x - y)| \leq L|x - y|$, for all $x, y \in [0, a]$ with $|x - x^*| \leq \tilde{\eta}$, $|y - x^*| \leq \tilde{\eta}$. Proposition 2.2 follows directly from Fact III with $M = \frac{k_2 + 1 - |1 - k_2|}{k_2}$ and $1 - g = k_2$. Indeed, by using induction and (2.11) we get:

$$V(x(t), w(t)) \leq \lambda^t V(x(0), w(0)), \text{ for all } t \geq 0$$

where $\lambda := \max\left(M^{-1} + L, |g| + ML|1 - g|\right) < 1$. ◁

**Proof of Proposition 2.1:** The proof of Proposition 2.1 is the same with the proof of Proposition 2.2, with the only difference that we use Fact III with $F = f'(x^*)$, $g = 1 - q$ and $M = \frac{|q| + 1 - |1 - q|}{|q|(1 - |1 - q|) + |k_2 - q^2|}$. ◁

## 3. Global Results

This section is devoted to the study of the (global) Input-to-State Stability (ISS) of system (1.4) with respect to the input $v \in \Re_+$ (uncontrollable inflow). More specifically, we study system (1.4) under the assumption:

$$x^* + v^* + b_{\max} < a \qquad (3.1)$$

We also assume that the uncertain function $f(d,x)$ satisfies the following assumption:

**(H2)** *There exist constants* $r \leq \dfrac{b_{\min} - u^*}{a - v^* - x^* - b_{\max}}$ *with* $0 < k_1 + k_2 + r < 2$, $\lambda_i \in [0,1)$, $\gamma_i \in (0, 1 - \lambda_i)$ ($i = 1,2$), $L \in [0,1)$, $q \in (0,1]$, $M \in (1, +\infty)$, *such that the following inequalities hold:*

$$-\frac{k_1 + LM^{-1}}{\beta + M^{-1}}(x - x^*) + \frac{\lambda b_{\min} + u^* - (1+\lambda)\max(b_{\min}, a - v^* - x)}{\beta + M^{-1}} \leq f(d,x) + x^* - a \leq 0$$
$$\text{for all } (d,x) \in D \times [a - v^* - b_{\max}, a] \qquad (3.2)$$

$$\left(f(d,x) + x^* - a\right)\left(\beta - M^{-1}\right) \leq \left(LM^{-1} - k_1\right)(x - x^*) + u^* - \lambda b_{\min} - (1-\lambda)b_{\max}$$
$$\text{for all } (d,x) \in D \times [a - v^* - b_{\max}, a] \qquad (3.3)$$

$$-(\lambda_1 q^{-1} + 1)(x - x^*) \leq P\!\left(u^* + r(x - x^*)\right) - f(d,x) + v^* \leq (\lambda_1 - 1)(x - x^*),$$
$$\text{for all } (d,x) \in D \times [x^*, a - v^* - b_{\min}] \qquad (3.4)$$

$$\left|(1-\beta)P\!\left(u^* + r(x - x^*)\right) + \beta f(d,x) - \beta v^* - (\beta - k_1)(x - x^*) - u^*\right| \leq \gamma_1(x - x^*),$$
$$\text{for all } (d,x) \in D \times [x^*, a - v^* - b_{\min}] \qquad (3.5)$$

$$(\lambda_2 - 1)(x - x^*) \leq P\!\left(u^* + r(x - x^*)\right) - f(d,x) + v^* \leq -(\lambda_2 q + 1)(x - x^*), \ \forall (d,x) \in D \times [0, x^*] \qquad (3.6)$$

$$\left|(1-\beta)P\!\left(u^* + r(x - x^*)\right) + \beta f(d,x) - \beta v^* - (\beta - k_1)(x - x^*) - u^*\right| \leq \gamma_2 q |x - x^*|, \ \forall (d,x) \in D \times [0, x^*] \qquad (3.7)$$

*where* $\beta := k_1 + k_2 + r$ *and* $P : \Re \to \Re$ *is defined in (2.2).*

Assumption (H2) is a set of sector-like conditions for the uncertain function $f(d,x)$. Figure 2 shows the allowable values for the uncertain function $f(d,x)$, as determined by assumption (H2) for $x^* = 10$, $u^* = 1$, $v^* = 2.678794$, $a = 16.8$, $\beta = q = 1$, $r = -0.98$, $k_1 = 0.9$, $k_2 = 1.08$, $L = 0.99$, $M = 1.025$, $\lambda_2 = 0.6$, $\gamma_2 = 0.39$, $\lambda_1 = 0.82$, $\gamma_1 = 0.17$, $b_{\min} = 0$ and $b_{\max} = 3.1$.

Assumption (H2) allows us to prove the following technical result.

**Lemma 3.1:** *Let* $a > 0$, $0 \leq b_{\min} < b_{\max}$ *be constants and let* $f \in C^0(D \times [0,a]; \Re)$ *be a function satisfying (1.2) for which assumptions (H1), (H2) and inequality (3.1) hold. Let* $V : [0,a] \times [b_{\min}, b_{\max}] \to \Re_+$ *be the function defined by*

$$V(y, w) := g(y - x^*) + M\left|w - P(u^* + r(y - x^*))\right| \qquad (3.8)$$

*where* $M > 1$ *is the constant involved in assumption (H2) and* $g : \Re \to \Re_+$ *is the function defined by* $g(x) := x$ *for* $x \geq 0$ *and* $g(x) := -qx$ *for* $x < 0$. *Then the following inequality holds for all* $(d, v, y, w) \in D \times \Re_+ \times [0,a] \times [b_{\min}, b_{\max}]$:

$$V\left(y-f(d,y)+\min(w+v,a-y),P\left(w+(k_1+k_2)(f(d,y)-\min(w+v,a-y))-k_2(y-x^*)\right)\right) \quad (3.9)$$
$$\leq \lambda V(y,w)+\gamma\left|v-v^*\right|$$

where $\lambda := \max\left(M^{-1}+|1-\beta|, L, \max_{i=1,2}(\lambda_i+M\gamma_i)\right)$ and $\gamma := 1+M|k_1+k_2|+M|r|$.

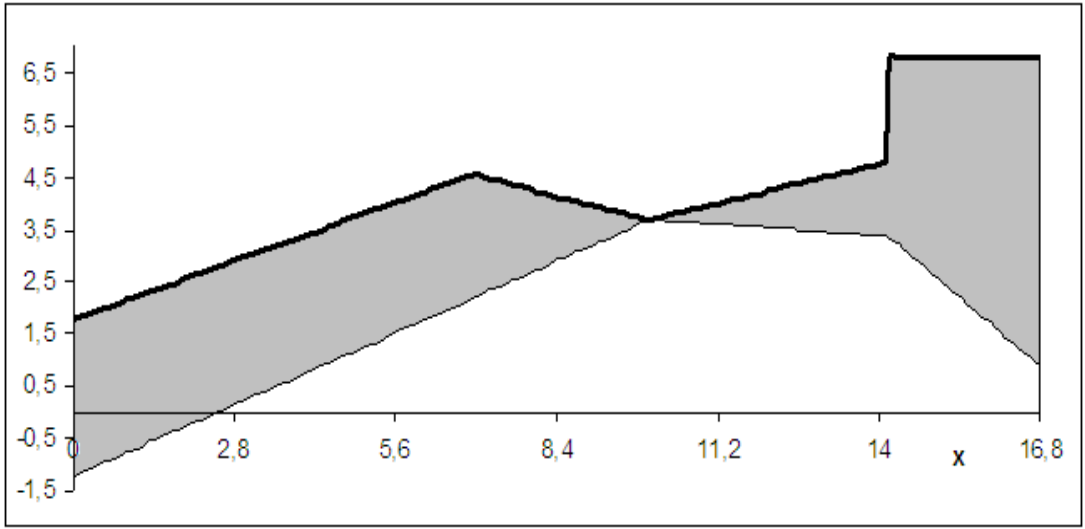

**Fig. 2:** The grey area shows the allowable values for the uncertain function $f(d,x)$, as determined by assumption (H2) for $x^* = 10$, $u^* = 1$, $v^* = 2.678794$, $a = 16.8$, $\beta = q = 1$, $r = -0.98$, $k_1 = 0.9$, $k_2 = 1.08$, $L = 0.99$, $M = 1.025$, $\lambda_2 = 0.6$, $\gamma_2 = 0.39$, $\lambda_1 = 0.82$, $\gamma_1 = 0.17$, $b_{\min} = 0$ and $b_{\max} = 3.1$.

**Proof:** Let arbitrary $(d,v,y,w) \in D \times \Re_+ \times [0,a] \times [b_{\min}, b_{\max}]$ and define:

$$V^+ := V\left(y-f(d,y)+\min(w+v,a-y), P\left(w+(k_1+k_2)(f(d,y)-\min(w+v,a-y))-k_2(y-x^*)\right)\right)$$
$$\tilde{V} := V\left(y-f(d,y)+\min(w+v^*,a-y), P\left(w+(k_1+k_2)(f(d,y)-\min(w+v^*,a-y))-k_2(y-x^*)\right)\right). \quad (3.10)$$

In what follows, we make extensive use of the inequalities:

$$|P(x)-P(y)| \leq |x-y|, \text{ for all } x,y \in \Re \quad (3.11)$$

$$|\min(u,a-y)-\min(w,a-y)| \leq |u-w|, \text{ for all } u,w,y \in \Re \quad (3.12)$$

$$|g(x)-g(y)| \leq |x-y|, \text{ for all } x,y \in \Re. \quad (3.13)$$

Definitions (3.8), (3.10) imply that:

$$V^+ = V\left(y-f(d,y)+\min(w+v,a-y), P\left(w+(k_1+k_2)(f(d,y)-\min(w+v,a-y))-k_2(y-x^*)\right)\right)$$
$$= g\left(y-f(d,y)+\min(w+v,a-y)-x^*\right) + \quad (3.14)$$
$$M\left|P\left(w+(k_1+k_2)(f(d,y)-\min(w+v,a-y))-k_2(y-x^*)\right)-P\left(u^*+r(y-f(d,y)+\min(w+v,a-y)-x^*)\right)\right|$$

$$\begin{aligned}
\tilde{V} &= V\big(y - f(d,y) + \min(w + v^*, a - y), P(w + (k_1 + k_2)(f(d,y) - \min(w + v^*, a - y)) - k_2(y - x^*))\big) \\
&= g\big(y - f(d,y) + \min(w + v^*, a - y) - x^*\big) + \\
&\quad M\big|P(w + (k_1 + k_2)(f(d,y) - \min(w + v^*, a - y)) - k_2(y - x^*)) - P(u^* + r(y - f(d,y) + \min(w + v^*, a - y) - x^*))\big|
\end{aligned}$$
(3.15)

Using (3.12) and (3.13) we obtain:

$$\begin{aligned}
&g\big(y - f(d,y) + \min(w + v, a - y) - x^*\big) \\
&\leq g\big(y - f(d,y) + \min(w + v^*, a - y) - x^*\big) + \big|\min(w + v, a - y) - \min(w + v^*, a - y)\big| \\
&\leq g\big(y - f(d,y) + \min(w + v^*, a - y) - x^*\big) + |v - v^*|.
\end{aligned}$$
(3.16)

Using (3.11), (3.12) and (3.13) we obtain:

$$\begin{aligned}
&\big|P(w + (k_1 + k_2)(f(d,y) - \min(w + v, a - y)) - k_2(y - x^*)) - P(u^* + r(y - f(d,y) + \min(w + v, a - y) - x^*))\big| \\
&\leq \big|P(w + (k_1 + k_2)(f(d,y) - \min(w + v^*, a - y)) - k_2(y - x^*)) - P(u^* + r(y - f(d,y) + \min(w + v^*, a - y) - x^*))\big| \\
&\quad + \big|P(w + (k_1 + k_2)(f(d,y) - \min(w + v, a - y)) - k_2(y - x^*)) - P(w + (k_1 + k_2)(f(d,y) - \min(w + v^*, a - y)) - k_2(y - x^*))\big| \\
&\quad + \big|P(u^* + r(y - f(d,y) + \min(w + v, a - y) - x^*)) - P(u^* + r(y - f(d,y) + \min(w + v^*, a - y) - x^*))\big| \\
&\leq \big|P(w + (k_1 + k_2)(f(d,y) - \min(w + v^*, a - y)) - k_2(y - x^*)) - P(u^* + r(y - f(d,y) + \min(w + v^*, a - y) - x^*))\big| \\
&\quad + (|k_1 + k_2| + |r|)\big|\min(w + v^*, a - y) - \min(w + v, a - y)\big| \\
&\leq \big|P(w + (k_1 + k_2)(f(d,y) - \min(w + v^*, a - y)) - k_2(y - x^*)) - P(u^* + r(y - f(d,y) + \min(w + v^*, a - y) - x^*))\big| \\
&\quad + (|k_1 + k_2| + |r|)|v - v^*|.
\end{aligned}$$
(3.17)

Combining (3.14), (3.15), (3.16) and (3.17) we get:

$$V^+ \leq \tilde{V} + \big(1 + M|k_1 + k_2| + M|r|\big)|v - v^*|$$
(3.18)

We next evaluate $\tilde{V}$. We distinguish three different cases.

CASE 1: $w + v^* \geq a - y$.

Consider first the case $w + v^* \geq a - y$, which necessarily implies $y \geq a - v^* - b_{\max}$. Since $\min(w + v^*, a - y) = a - y$, we obtain from (3.15) (using (3.11)):

$$\tilde{V} \leq g\big(a - f(d,y) - x^*\big) + M\big|w + (k_1 + k_2 + r)(f(d,y) + y - a) - (k_2 + r)(y - x^*) - u^*\big|.$$
(3.19)

Since $r \leq \dfrac{b_{\min} - u^*}{a - v^* - x^* - b_{\max}}$ and $y \geq a - v^* - b_{\max}$, we get $u^* + r(y - x^*) \leq b_{\min}$. By virtue of (3.1) and the fact that $y \geq a - v^* - b_{\max}$, we conclude that $y \geq x^*$ and consequently $g(y - x^*) = |y - x^*|$. The two previous observations in conjunction with definitions (2.2), (3.8) imply that:

$$V(y, w) = y - x^* + M(w - b_{\min}).$$
(3.20)

The right hand side inequality (3.2) directly implies that

$$g(a - f(d,y) - x^*) = a - f(d,y) - x^*. \tag{3.21}$$

By virtue of (3.20) and (3.21), the inequality

$$g(a - f(d,y) - x^*) + M\left|w + (k_1 + k_2 + r)(f(d,y) + y - a) - (k_2 + r)(y - x^*) - u^*\right| \le LV(y,w) \tag{3.22}$$

is equivalent to the two following inequalities:

$$w + \beta(f(d,y) + y - a) - (\beta - k_1)(y - x^*) - u^* \le LM^{-1}(y - x^*) + L(w - b_{\min}) - M^{-1}(a - f(d,y) - x^*)$$

and

$$w + \beta(f(d,y) + y - a) - (\beta - k_1)(y - x^*) - u^* \ge -LM^{-1}(y - x^*) - L(w - b_{\min}) + M^{-1}(a - f(d,y) - x^*).$$

Since $\max(b_{\min}, a - v^* - y) \le w \le b_{\max}$ we conclude that the two above inequalities hold, provided that the two following inequalities hold:

$$(1 - L)b_{\max} + \beta(f(d,y) + y - a) - u^* \le (\beta - k_1 + LM^{-1})(y - x^*) - Lb_{\min} - M^{-1}(a - f(d,y) - x^*) \tag{3.23a}$$

and

$$(1 + L)\max(b_{\min}, a - v^* - y) + \beta(f(d,y) + y - a) - u^* \ge (\beta - k_1 - LM^{-1})(y - x^*) + Lb_{\min} + M^{-1}(a - f(d,y) - x^*) \tag{3.23b}$$

The fact that inequalities (3.23) hold, is a direct consequence of inequalities (3.2) and (3.3). Thus inequality (3.22) holds. Combining inequalities (3.19) and (3.22), we obtain:

$$\tilde{V} \le LV(y,w). \tag{3.24}$$

<u>CASE 2:</u> $w + v^* \le a - y$ and $y \ge x^*$.

Inequality $w + v^* \le a - y$ necessarily implies $y \le a - v^* - b_{\min}$. Since $\min(w + v^*, a - y) = w + v^*$, we obtain from (3.15) (using (3.11)):

$$\tilde{V} \le g(y - f(d,y) + w + v^* - x^*) + M\left|(1 - \beta)w + \beta f(d,y) - \beta v^* - (\beta - k_1)(y - x^*) - u^*\right|. \tag{3.25}$$

Inequality (3.13) in conjunction with (3.25) implies the following inequality:

$$\tilde{V} \le g(y - x^* - f(d,y) + v^* + P(u^* + r(y - x^*)))$$
$$+ M\left|(1 - \beta)P(u^* + r(y - x^*)) + \beta f(d,y) - \beta v^* - (\beta - k_1)(y - x^*) - u^*\right| \tag{3.26}$$
$$+ (1 + M|1 - \beta|)|w - P(u^* + r(y - x^*))|.$$

Using inequality (3.5) in conjunction with inequality (3.26), we obtain:

$$\tilde{V} \le g(y - x^* - f(d,y) + v^* + P(u^* + r(y - x^*))) + M\gamma_1|y - x^*|$$
$$+ (1 + M|1 - \beta|)|w - P(u^* + r(y - x^*))|. \tag{3.27}$$

Inequality (3.4) in conjunction with the fact that $y \ge x^*$, directly implies that

$$g(y-x^* - f(d,y) + v^* + P(u^* + r(y-x^*))) \leq \lambda_1 |y-x^*|. \tag{3.28}$$

Using inequalities (3.27), (3.28) in conjunction with the fact that $y \geq x^*$ (which implies that $g(y-x^*) = |y-x^*|$) and definition (3.8), we get:

$$\widetilde{V} \leq \max(\lambda_1 + M\gamma_1, M^{-1} + |1-\beta|) V(y,w). \tag{3.29}$$

<u>CASE 3:</u> $w + v^* \leq a - y$ and $y \leq x^*$.

Since $\min(w+v^*, a-y) = w+v^*$, we obtain (3.25) from (3.15) (using (3.11)). Inequality (3.13) in conjunction with (3.25) implies inequality (3.26). Using inequality (3.7) in conjunction with inequality (3.26), we obtain:

$$\begin{aligned}\widetilde{V} &\leq g(y-x^* - f(d,y) + v^* + P(u^* + r(y-x^*))) + M\gamma_2 q|y-x^*| \\ &+ (1+M|1-\beta|)|w - P(u^* + r(y-x^*))|.\end{aligned} \tag{3.30}$$

Inequality (3.6) in conjunction with the fact that $y \leq x^*$, directly implies that

$$g(y-x^* - f(d,y) + v^* + P(u^* + r(y-x^*))) \leq \lambda_2 g(y-x^*). \tag{3.31}$$

Using inequalities (3.30), (3.31) in conjunction with the fact that $y \leq x^*$ (which implies that $g(y-x^*) = q|y-x^*|$) and definition (3.8), we obtain:

$$\widetilde{V} \leq \max(\lambda_2 + M\gamma_2, M^{-1} + |1-\beta|) V(y,w). \tag{3.32}$$

Combining inequalities (3.18), (3.24), (3.29) and (3.32), we obtain inequality (3.9) with $\lambda := \max\left(M^{-1} + |1-\beta|, L, \max_{i=0,1,2}(\lambda_i + M\gamma_i)\right)$ and $\gamma := 1 + M|k_1 + k_2| + M|r|$. The proof is complete. ◁

The following result provides sufficient conditions for the ISS property with respect to the external input $v \in \Re_+$ uniformly in $d \in D$ for the closed-loop system (1.4). Notice that the gain of the external input $v \in \Re_+$ is linear and is explicitly given. Moreover, for $v \equiv v^*$ we have exponential convergence with rate which is explicitly estimated.

**Theorem 3.2:** *Let $a > 0$, $0 \leq b_{\min} < b_{\max}$ be constants and let $f \in C^0(D \times [0,a]; \Re)$ be a function satisfying (1.2) for which assumptions (H1), (H2) and inequality (3.1) hold. Consider system (1.4) and suppose that*

$$1 + \min_{i=1,2}\left(\frac{1-\lambda_i}{\gamma_i}\right)|1-\beta| < \min_{i=1,2}\left(\frac{1-\lambda_i}{\gamma_i}\right), \quad \frac{1}{1-|1-\beta|} < M < \min_{i=1,2}\left(\frac{1-\lambda_i}{\gamma_i}\right) \tag{3.33}$$

*where $\beta := k_1 + k_2 + r$. Then the following estimate holds for the solution of (1.4) corresponding to arbitrary inputs $\{d(i) \in D\}_{i=0}^{\infty}$ and $\{v(i) \in \Re_+\}_{i=0}^{\infty}$:*

$$\begin{aligned}\max(q(x^* - x(t)), x(t) - x^*) &\leq \\ \lambda^t M(1+|r|+|k_1|+|k_2|)\left(|x(0)-x^*| + |y(0)-x^*| + |w(0)-u^*|\right) &+ \frac{\gamma}{1-\lambda}\max_{i=0,...,t-1}\left(|v(i)-v^*|\right),\end{aligned} \quad \text{for all } t \geq 1 \tag{3.34}$$

*where $\lambda := \max\left(M^{-1} + |1-\beta|, L, \max_{i=1,2}(\lambda_i + M\gamma_i)\right)$ and $\gamma := 1 + M|k_1 + k_2| + M|r|$.*

**Proof:** Consider the solution of (1.4) corresponding to arbitrary inputs $\{d(i) \in D\}_{i=0}^{\infty}$ and $\{v(i) \in \Re_+\}_{i=0}^{\infty}$. It follows that the following equations hold for all $t \geq 1$:

$$y(t) = x(t-1) \quad , \quad w(t) = u(t-1) \tag{3.35}$$

$$y(t+1) = y(t) - f(d(t-1), y(t)) + \min(u(t-1) + v(t-1), a - y(t)) \tag{3.36}$$

$$w(t+1) = P\big(w(t) + (k_1 + k_2)(f(d(t-1), y(t)) - \min(w(t) + v(t-1), a - y(t))) - k_2(y(t) - x^*)\big). \tag{3.37}$$

Using the results of Lemma 3.1, we conclude that the following inequality holds for $t \geq 1$:

$$V(y(t+1), w(t+1)) \leq \lambda V(y(t), w(t)) + \gamma |v(t-1) - v^*| \tag{3.38}$$

with $\lambda := \max\left(M^{-1} + |1-\beta|, L, \max_{i=1,2}(\lambda_i + M\gamma_i)\right)$ and $\gamma := 1 + M|k_1 + k_2| + M|r|$. Inequalities (3.33) imply that

$$\lambda = \max\left(M^{-1} + |1-\beta|, L, \max_{i=1,2}(\lambda_i + M\gamma_i)\right) < 1. \tag{3.39}$$

Using induction and inequality (3.38), we conclude that the following inequality holds for all $N \geq 1$:

$$V(y(N+1), w(N+1)) \leq \lambda^N V(y(1), w(1)) + \gamma \sum_{i=0}^{N-1} \lambda^{N-1-i} |v(i) - v^*|. \tag{3.40}$$

Equations (3.35) in conjunction with definition (3.8) and inequality (3.40) imply the following estimate for all $N \geq 1$:

$$g(x(N) - x^*) \leq \lambda^N V(x(0), u(0)) + \frac{\gamma}{1-\lambda} \max_{i=0,\ldots,N-1}\left(|v(i) - v^*|\right). \tag{3.41}$$

Inequality (3.34) is a direct consequence of estimate (3.41) and definition (3.8). The proof is complete. ◁

**Remark 3.3:** Theorem 3.2 guarantees the ISS property for system (1.4) with respect to the input $v \in \Re_+$ uniformly in $d \in D$. However, since we are most interested in the component $x$ of the solution of (1.4) we have provided only estimate (3.34). Similar Sontag-like estimates hold for all components of the solution of (1.4).

**Remark 3.4:** Saturation of traffic systems is a very important phenomenon, which must be avoided. It is reasonable to adopt the following definition for the saturation phenomenon of the traffic system (1.1) under the effect of the PI-regulator (1.3):

> "We say that the traffic network becomes saturated for inputs $\{d(i) \in D\}_{i=0}^{\infty}$, $\{v(i) \in \Re_+\}_{i=0}^{\infty}$ and initial condition $(x(0), y(0), w(0)) \in [0,a]^2 \times \Re$ if for every $N \geq 1$ there exists $t \geq N$ such that the corresponding solution of (1.4) satisfies $x(t) = a$."

Estimate (3.34) guarantees that:

> "for every pair of inputs $\{d(i) \in D\}_{i=0}^{\infty}$, $\{v(i) \in \Re_+\}_{i=0}^{\infty}$ with $\sup_{i \geq 0}\left(|v(i) - v^*|\right) < \frac{1-\lambda}{\gamma} a$ and for every initial condition $(x(0), y(0), w(0)) \in [0,a]^2 \times \Re$, the traffic network cannot become saturated, provided that the assumptions of Theorem 3.2 hold."

Thus, Theorem 3.2 allows us to estimate the range of values for the uncontrollable inflow $\{v(i) \in \Re_+\}_{i=0}^{\infty}$, for which the traffic network cannot become saturated for any initial condition and for any uncertainty $\{d(i) \in D\}_{i=0}^{\infty}$.

Finally, we end this section by providing a set of necessary conditions for the input-to-state stabilizability by means of the PI regulator.

**Theorem 3.5:** *Let $a > 0$, $0 \leq b_{min} < b_{max}$ be constants and let $f \in C^0(D \times [0,a]; \Re)$ be a function satisfying (1.2). Suppose that (3.1) holds and that there exist constants $M > 0$, $\lambda \in (0,1)$ and $\gamma \geq 0$, such that inequality (3.34) holds for the solution of (1.4) corresponding to arbitrary inputs $\{d(i) \in D\}_{i=0}^{\infty}$ and $\{v(i) \in \Re_+\}_{i=0}^{\infty}$. Then $k_2 > 0$ and the following properties hold for every $d \in D$:*

"either the equation $f(d,y) = \min(b_{max} + v^*, a - y)$ has no solution $y \in [0,a]$

or all solutions $y \in [0,a]$ of the equation $f(d,y) = \min(b_{max} + v^*, a - y)$ are in $(x^*, a]$" (3.42)

"all solutions $y \in [0,a]$ of the equation $f(d,y) = \min(b_{min} + v^*, a - y)$ are in $[0, x^*)$." (3.43)

**Proof:** As in the proof of Theorem 3.2, we consider the solution of (1.4) with $v \equiv v^*$ corresponding to constant input $d(i) \equiv d \in D$. It follows that the following equations hold for all $t \geq 1$:

$$y(t) = x(t-1) \quad , \quad w(t) = u(t-1) \tag{3.44}$$

$$y(t+1) = y(t) - f(d, y(t)) + \min(w(t) + v^*, a - y(t)) \tag{3.45}$$

$$w(t+1) = P\Big(w(t) + (k_1 + k_2)(f(d, y(t)) - \min(w(t) + v^*, a - y(t))) - k_2(y(t) - x^*)\Big) \tag{3.46}$$

Suppose that the system of equations
$$f(d,y) = \min(w + v^*, a - y) \; ; \; w = P\big(w - k_2(y - x^*)\big) \tag{3.47}$$

admits a solution $(y, w) = (y^*, w^*) \in [0,a] \times [b_{min}, b_{max}]$ other than $(x^*, u^*)$. Then the solution of (1.4) with $v \equiv v^*$ corresponding to constant input $d(i) \equiv d \in D$ starting from any initial condition $(x(0), y(0), w(0)) \in [0,a]^2 \times \Re$ with $x(0) = y^*$, $w^* = \max\big(b_{min}, \min\big(w(0) + k_1 y(0) - (k_1 + k_2) y^* + k_2 x^*, b_{max}\big)\big)$ would satisfy $x(t) = y^*$ for all $t \geq 0$, which contradicts (3.34) (which implies $\lim x(t) = x^*$). However, notice that, for each fixed $y(0) \in [0,a]$, there exists at least one $w(0) \in \Re$ for which $w^* = \max\big(b_{min}, \min\big(w(0) + k_1 y(0) - (k_1 + k_2) y^* + k_2 x^*, b_{max}\big)\big)$. Therefore, the system of equations (3.47) must have a unique solution, which necessarily is the equilibrium point $(y, w) = (x^*, u^*) \in [0,a] \times [b_{min}, b_{max}]$.

Consequently, we cannot have a solution $(y, w) \in [0,a] \times [b_{min}, b_{max}]$ of the system of equations (3.47) for which $w = b_{min}$. This means that either the equation $f(d, y) = \min(b_{min} + v^*, a - y)$ has no solution $y \in [0, a]$ or that every solution $y \in [0,a]$ of the equation $f(d,y) = \min(b_{min} + v^*, a - y)$ satisfies $k_2(y - x^*) < 0$. Clearly, the latter means that $k_2 \neq 0$ and that all solutions $y \in [0,a]$ of the equation $f(d,y) = \min(b_{min} + v^*, a - y)$ must be either in $[0, x^*)$ or in $(x^*, a]$. However, there exists at least one solution $y \in (0, x^*)$ of the equation $f(d,y) = \min(b_{min} + v^*, a - y)$: the existence of such a solution follows from the inequalities $f(d, x^*) = u^* + v^* > b_{min} + v^*$ and $f(d, 0) = 0 < b_{min} + v^*$, which imply that there exists $y \in (0, x^*)$ for which $f(d,y) = b_{min} + v^*$ (and since (3.1) holds, we get $v^* + b_{min} \leq a - x^* < a - y$ and consequently $y \in (0, x^*)$ is a solution of $f(d,y) = \min(b_{min} + v^*, a - y)$). Therefore, we conclude that all solutions $y \in [0,a]$ of the equation $f(d,y) = \min(b_{min} + v^*, a - y)$ must be in $[0, x^*)$ and that $k_2 > 0$.

Similarly, we cannot have a solution $(y, w) \in [0,a] \times [b_{min}, b_{max}]$ of the system of equations (3.47) for which $w = b_{max}$. This means that either the equation $f(d,y) = \min(b_{max} + v^*, a - y)$ has no solution $y \in [0,a]$ or that every solution $y \in [0,a]$ of the equation $f(d,y) = \min(b_{max} + v^*, a - y)$ satisfies $k_2(y - x^*) > 0$. Since $k_2 > 0$, the latter means that all solutions $y \in [0,a]$ of the equation $f(d,y) = \min(b_{max} + v^*, a - y)$ must be in $(x^*, a]$.

The proof is complete. ◁

**Remark 3.6:** The existence of a solution $y^* \in (x^*, a]$ of the equation $f(d, y) = \min(b_{\min} + v^*, a - y)$ for certain $d \in D$ implies the existence of an equilibrium point for system (1.4), which cannot be removed by the control action of the PI-regulator, no matter what the values of $k_1, k_2$ are (namely the equilibrium point $(y^*, y^*, b_{\min})$). The reader should notice that the existence of a solution $y^* \in (x^*, a]$ of the equation $f(d, y) = \min(b_{\min} + v^*, a - y)$ depends only on the properties of the outflow function $f(d, y)$ and the parameters $v^*, x^*, a, b_{\min}$. Therefore, since the control practitioner cannot change the properties of the outflow function $f(d, y)$ and the parameters $v^*, x^*, a, b_{\min}$ (they are characteristic of the given traffic system), global stabilization (in the sense of [6]) cannot be achieved in this case, no matter what the values of $k_1, k_2$ are.

## 4. Illustrative Examples

All examples in this section assume that the outflow function $f(d, x)$ is not uncertain (i.e., $f(d, x) = f(x)$) and is of the form:

$$f(x) = px \exp(-cx^\delta)$$

where $p \in (0,1]$, $c, \delta > 0$ are constants. This is the most usual form for the outflow function $f(d, x)$ for traffic systems.

**Example 4.1:** Consider system (1.4) with $f(x) = x \exp(-x/10)$, $x^* = 10$, $u^* = 1$, $v^* = 2.678794$, $a = 16.8$, $k_1 = 0.9$, $k_2 = 1.08$, $b_{\min} = 0$ and $b_{\max} = 3.1$. The regions of attractions of the corresponding system (2.1) as predicted by Propositions 2.1 and 2.2 are shown in Figure 3. The size of the resulting region of attraction may be not sufficiently big for particular application; in this case, the estimated region of attraction may provide the necessary basis for a numerical elaboration of the exact region of attraction, see [2,5,15] for details.

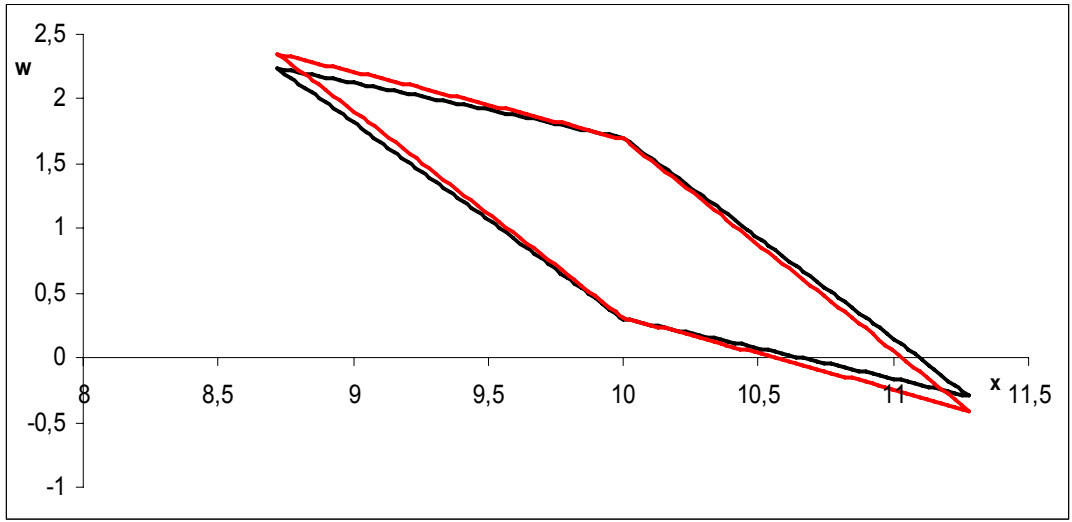

**Fig. 3:** The red line and the black line show the boundaries of the regions of attraction for system (2.1) with

$f(x) = x \exp\left(-\dfrac{x}{10}\right)$, $x^* = 10$, $u^* = 1$, $v^* = 2.678794$, $a = 16.8$, $k_1 = 0.9$, $k_2 = 1.08$, $b_{\min} = 0$ and

$b_{\max} = 3.1$, as predicted by Proposition 2.2 and Proposition 2.1, respectively.

On the other hand, we can use Theorem 3.2 and check whether we have global exponential stability. Indeed, the allowable values for the outflow function $f(x)$, as determined by assumption (H2) for $x^* = 10$, $u^* = 1$, $v^* = 2.678794$, $a = 16.8$, $\beta = q = 1$, $r = -0.98$, $k_1 = 0.9$, $k_2 = 1.08$, $L = 0.99$, $M = 1.025$, $\lambda_2 = 0.6$, $\gamma_2 = 0.39$, $\lambda_1 = 0.82$, $\gamma_1 = 0.17$, $b_{\min} = 0$ and $b_{\max} = 3.1$ was shown in Figure 2. Figure 4 shows that in this case the graph of $f(x) = x \exp(-x/10)$ lies completely in the area of the allowable values for the outflow function $f(x)$. Therefore,

Theorem 3.2 guarantees global exponential stability of system (1.4) with $f(x) = x\exp(-x/10)$, $x^* = 10$, $u^* = 1$, $v \equiv v^* = 2.678794$, $a = 16.8$, $k_1 = 0.9$, $k_2 = 1.08$, $b_{\min} = 0$ and $b_{\max} = 3.1$.

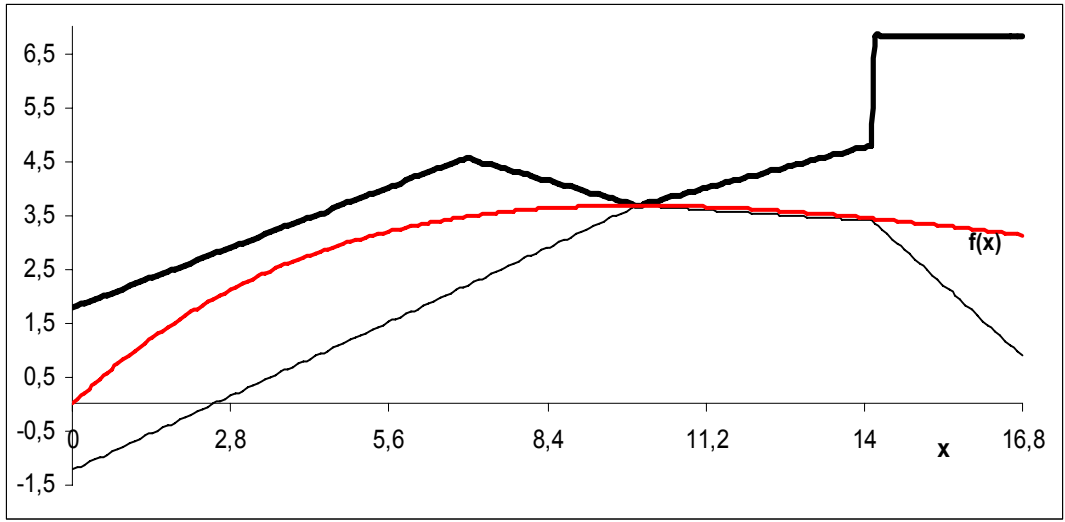

**Fig. 4:** The graph of $f(x) = x\exp(-x/10)$ is shown in red color.

In this example, the estimation of the region of attraction by means of Propositions 2.1 and 2.2 is very conservative. The reason of the conservatism lies in the saturation effects of the input and the inflow: in this example the saturation effects of the input and the inflow are strongly stabilizing. On the other hand, Propositions 2.1 and 2.2 are based on estimates for which no saturation effects of the input and the inflow are present and, consequently, cannot take advantage of the nonlinear stabilizing nature of the saturation effects. ◁

**Example 4.2:** Consider system (1.4) with $f(x) = x\exp(-x^2/10)$, $x^* = 3$, $u^* = 0.219709$, $v^* = 1$, $a = 20$, $b_{\min} = 0$ and $b_{\max} = 3$. For this system one cannot hope for global stabilization by means of the PI-regulator, since there exists a solution $y^* \in (3.5, 3.6)$ of the equation $f(y) = \min(1, 20 - y)$ (recall Remark 3.6). However, we can apply the PI-regulator for local stabilization. Indeed, the selection $k_1 = k_2 = 1$ guarantees local exponential stabilization. The regions of attractions of the corresponding system (2.1) as predicted by Propositions 2.1, 2.2 are shown in Figure 5.

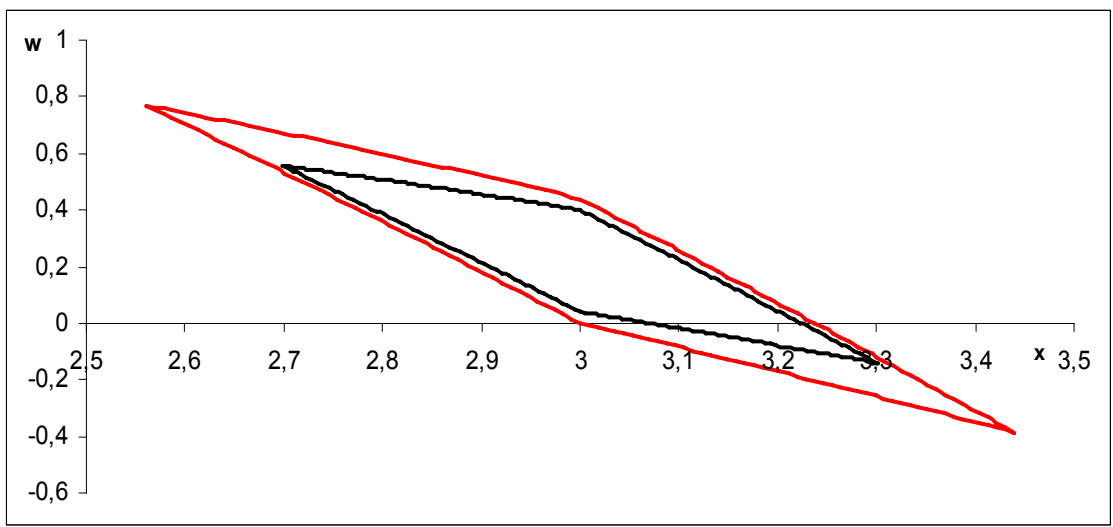

**Fig. 5:** The red line and the black line show the boundaries of the regions of attraction for system (2.1) with $f(x) = x\exp(-x^2/10)$, $x^* = 3$, $u^* = 0.219709$, $v^* = 1$, $a = 20$, $k_1 = k_2 = 1$, $b_{\min} = 0$ and $b_{\max} = 3$ as predicted by Proposition 2.2 and Proposition 2.1, respectively.

In this case, the estimate of the region of attraction for system (2.1) provided by Proposition 2.2 is not conservative: the existence of a solution $y^* \in (3.5, 3.6)$ of the equation $f(y) = \min(1, 20 - y)$ (and consequently the existence of the equilibrium point $(y^*, 0)$ for system (2.1)) implies that the constant $\rho > 0$ involved in the estimate of the region of attraction $\Omega_\rho = \left\{ (x, w) \in [0, a] \times \Re : \left| x - x^* \right| + \frac{k_2 + 1 - |1 - k_2|}{k_2} \left| w + v^* - f(x) + k_2(x - x^*) \right| < \rho \right\}$ cannot exceed the number $0.6$, while the constant $\rho > 0$ computed by (2.7) is $0.439418$.  ◁

## 5. Concluding Remarks

This work has provided necessary conditions and sufficient conditions for the (global) Input-to-State Stability property of simple uncertain vehicular-traffic networks under the effect of a PI-regulator controller (Theorem 3.2 and Theorem 3.5). We have also studied the local stability properties for vehicular-traffic networks under the effect of PI-regulator control: the region of attraction of a locally exponentially stable equilibrium point was estimated by means of Lyapunov functions (Proposition 2.1 and Proposition 2.2). The obtained results were illustrated by means of two simple examples.

More remains to be done. One research direction is the application of PI-regulator control to traffic systems with uncertain input delays. Another research direction is the application of nonlinear feedback stabilizers to traffic systems. Both research directions will be the topic of future works.

**Acknowledgments**

The research leading to these results has received funding from the European Research Council under the European Union's Seventh Framework Programme (FP/2007-2013) / ERC Grant Agreement n. [321132].